\documentclass[10pt]{article}
\usepackage{graphicx}
\usepackage{tikz}
\usepackage[margin=1in]{geometry}
\usepackage{amsmath, amssymb, amsfonts, amsthm}
\newtheorem{lem}{\noindent {\bf Lemma}}[section]

\newtheorem{thm}{\noindent {\bf Theorem}}[section]

\newcounter{remark}
\newenvironment{remark}{\smallskip\noindent {\bf Remark \arabic{section}.\arabic{remark}.}}
{\addtocounter{remark}{1}\par}
\newcounter{example}
\catcode`@=11 
\date{}
\newcounter{defi}
\newenvironment{defi}{
\smallskip \noindent
{\bf
  Definition \arabic{section}.\arabic{defi}.
}}{\addtocounter{defi}{1}\par}
\newcommand{\ZZ}{\mathbb Z}
\newcommand{\NN}{\mathbb N}

\title{\bf Asymptotic property C of the wreath product $\ZZ \wr\ZZ$}
\author{\large Jingming Zhu$^\ast$\qquad  Yan Wu$^{\ast\ast}$
\footnote{
College of Mathematics Physics and Information Engineering, Jiaxing University, Jiaxing , 314001, P.R.China.
$^\ast$ E-mail: 122411741@qq.com $^\ast\ast$ E-mail:  yanwu@mail.zjxu.edu.cn}}
\date{}

\begin{document}
\maketitle
\begin{center}
\begin{minipage}{0.9\textwidth}
\noindent{\bf Abstract.}
Using the relationship between transfinite asymptotic dimension and asymptotic property C, we obtain that the wreath product $\ZZ \wr\ZZ$ has asymptotic property C. Specifically, we prove that the transfinite asymptotic dimension of the wreath product $\ZZ \wr\ZZ$ is no more than $\omega+1$.


\end{minipage}
\end{center}
\footnote{
This research was supported by
the National Natural Science Foundation of China under Grant (No.12071183,11871342)

}
\begin{section}{Introduction}\

In coarse geometry, asymptotic dimension of a metric space is an important concept which was defined by Gromov for studying asymptotic invariants of discrete groups \cite{Gromov}. This dimension can be considered as an asymptotic analogue of the Lebesgue covering dimension. As a large scale analogue of W.E. Haver¡¯s property C in dimension theory, A. Dranishnikov introduced the notion of asymptotic property C and proved that every metric
space of bounded geometry with asymptotic property C has property A~\cite{Dra00}.
E. Guentner,
R. Tessera and G. Yu \cite{Yu2012} introduced the notion of finite decomposition complexity to study topological
rigidity of manifolds, and proved that every metric space of bounded geometry with finite decomposition
complexity has property A \cite{Yu2013}.

 It is well known that every metric space with finite asymptotic dimension has asymptotic property C and finite decomposition complexity(\cite{Dra00}, \cite{Yu2013}). But the inverse is not true, which means that there exists some infinite asymptotic dimension metric space $X$ with
 asymptotic property C and finite decomposition complexity. Therefore how to classify the metric spaces with infinite asymptotic dimension into smaller categories becomes an interesting problem.  T. Radul defined trasfinite asymptotic dimension (trasdim) which can be viewed as transfinite extension for asymptotic dimension and proved that for every metric space $X$, $X$ has asymptotic property C is equivalent to trasdim$(X)\leq \alpha$ for some countable ordinal number $\alpha$ (see \cite{Radul2010}).

The relation between asymptotic property C and finite decomposition complexity was studied by A. Dranishnikov and M. Zarichnyi \cite{Mari2014}. There is no group of examples known which make a difference between asymptotic property C and finite decomposition complexity. The wreath product $\ZZ \wr\ZZ$  has finite decomposition complexity(see \cite{Yu2013},\cite{Yan2011}), but $\ZZ \wr\ZZ$  has infinite asymptotic dimension. So we are interested in the question whether $\ZZ \wr\ZZ$ has asymptotic property C.
In this paper, we prove that transfinite asymptotic dimension of the wreath product $\ZZ \wr\ZZ$ is no more than $\omega+1$. Consequently, $\ZZ \wr\ZZ$ has asymptotic property C. $\ZZ \wr\ZZ$ is the first finitely generated group we found with asymptotic property C and  infinite asymptotic dimension.

The paper is organized as follows: In Section 2, we recall some definitions and properties of transfinite asymptotic dimension, asymptotic property C
and wreath product. In
Section 3, we prove that transfinite asymptotic dimension of the wreath product $\ZZ \wr\ZZ$ is no more than $\omega+1$.
\end{section}

\begin{section}{Preliminaries}\

Our terminology concerning the asymptotic dimension follows from \cite{Bell2011} and for undefined terminology we refer to \cite{Radul2010}.
\begin{subsection}{Asymptotic Property C}

 Let~$(X, d)$ be a metric space and $U,V\subseteq X$. Let
\[
\text{diam}~ U=\text{sup}\{d(x,y): x,y\in U\}
\text{   and   }
d(U,V)=\text{inf}\{d(x,y): x\in U,y\in V\}.
\]
Let $R>0$ and $\mathcal{U}$ be a family of subsets of $X$, $\mathcal{U}$ is said to be \emph{$R$-bounded} if
\[
\text{diam}~\mathcal{U}\stackrel{\bigtriangleup}{=}\text{sup}\{\text{diam}~ U: U\in \mathcal{U}\}\leq R.
\]
$\mathcal{U}$ is said to be \emph{uniformly bounded} if there exists $R>0$
such that $\mathcal{U}$ is $R$-bounded.\\
Let $r>0$, $\mathcal{U}$ is said to be\emph{ $r$-disjoint} if
\[
d(U,V)\geq r~~~~~\text{for every}~ U,V\in \mathcal{U}\text{~and~}U\neq V.
\]
In this paper, we denote
$\bigcup\{U~|~U\in\mathcal{U}\}$ by $\bigcup\mathcal{U}$, denote $\{U~|~U\in\mathcal{U}_{1}\text{~or~}~U\in\mathcal{U}_{2}\}$
by $\mathcal{U}_{1}\cup\mathcal{U}_{2}$.

A metric space $X$ is said to have \emph{finite asymptotic dimension} if
there is an $n\in\NN$, such that for every $r>0$,
there exists a sequence of uniformly bounded families
$\{\mathcal{U}_{i}\}_{i=0}^{n}$ of subsets of $X$
such that the family
$\bigcup_{i=0}^{n}\mathcal{U}_{i}$ covers $X$ and each $\mathcal{U}_{i}$
is $r$-disjoint for $i=0,1,\cdots,n$. In this case, we say that asdim$(X)\leq n$.

A metric space $X$ is said to have \emph{asymptotic property C} if for every sequence $R_0<R_1<...$ of positive real numbers, there exist an $n\in\mathbb{N}$ and uniformly bounded families $\mathcal{U}_0,...,\mathcal{U}_n$ of subsets of $X$ such that each $\mathcal{U}_i$ is $R_i$-disjoint for $i=0,1,\cdots,n$ and the family $\bigcup_{i=0}^n\mathcal{U}_i$ covers $X$.\

\end{subsection}
\begin{subsection}{Transfinite Asymptotic Dimension}

In \cite{Radul2010}, T. Radul generalized asymptotic dimension of a metric space $X$ to transfinite asymptotic dimension which is denoted by trasdim$(X)$.

Let $Fin~\mathbb{N}$ denote the collection of all finite, nonempty
subsets of $\mathbb{N}$ and let $ M \subset Fin~\mathbb{N}$. For $\sigma\in \{\varnothing\}\bigcup Fin~\mathbb{N}$, let
$$M^{\sigma} = \{\tau\in Fin\mathbb{N} ~|~ \tau \cup \sigma \in M \text{ and } \tau \cap \sigma = \varnothing\}.$$

Let $M^a$ abbreviate $M^{\{a\}}$ for $a \in \NN$. Define ordinal number Ord$M$ inductively as follows:
\begin{eqnarray*}
\text{Ord}M = 0 &\Leftrightarrow& M = \varnothing,\\
\text{Ord}M \leq \alpha &\Leftrightarrow& \forall~ a\in \mathbb{N}, ~\text{Ord}M^a < \alpha,\\
\text{Ord}M = \alpha &\Leftrightarrow& \text{Ord}M \leq \alpha \text{ and } \text{Ord}M < \alpha \text{ is not true},\\
\text{Ord}M = \infty &\Leftrightarrow& \text{Ord}M \leq\alpha \text{ is not true for every ordinal number } \alpha.
\end{eqnarray*}

Given a metric space $(X, d)$, define the following collection:
\[
\begin{split}
A(X, d) = \{\sigma \in Fin\mathbb{N} |~&\text{ there are no uniformly bounded families } \mathcal{U}_i  \text{ for } i \in \sigma
 \\& \text{ such that each } \mathcal{U}_i
\text{ is } i\text{-disjoint and }\bigcup_{i\in\sigma}\mathcal{U}_i \text{~covers~} X\}.
\end{split}\]

The \emph{transfinite asymptotic dimension} of $X$ is defined as trasdim$(X)$=Ord$A(X, d)$.

\begin{remark}
It is not difficult to see that transfinite asymptotic dimension is a generalization of finite asymptotic dimension.
That is, trasdim $(X)\leq n$ if and only if asdim $(X)\leq n$ for each $n\in\NN$.

\end{remark}

\begin{lem}\rm(see \cite{Radul2010})
\label{lem:E}
Let $X$ be a metric space, $X$ has asymptotic property $C$ if and only if trasdim$(X)\leq\alpha$ for some countable ordinal number $\alpha$.
\end{lem}

\begin{lem}\rm(see \cite{Yan2018} )
\label{lem:trasdim}
Given a metric space $X$ with asdim$(X) = \infty$ and $k\in \mathbb{N}$, the following are equivalent:
\begin{itemize}
\item\rm trasdim$(X) \leq \omega+k$;
\item\rm For every $n\in\NN$, there exists $m(n)\in\NN,$ such that for every~ $d>0$, there are uniformly bounded families $\mathcal{U}_{-k},\mathcal{U}_{-k+1},\cdots,\mathcal{U}_{m(n)}$ such that $\mathcal{U}_i$ is $n$-disjoint for $i=-k,\cdots, 0$, $\mathcal{U}_j$ is $d$-disjoint for $j=1,2,\cdots, m(n)$ and
$\bigcup_{i=-k}^{m(n)}\mathcal{U}_i$ covers $X$. Moreover, $m(n)\rightarrow \infty$ as $n\rightarrow \infty$.
\end{itemize}
\end{lem}
\end{subsection}
\begin{subsection}{Wreath Product}

Let $S$ be a finite generating set for a group $G$, for any $ g\in G$, let $|g|_{S}$ to be
the length of the shortest word representing $g$ in elements of $S\cup S^{-1}$.
We say that $|\cdot|_{S}$ is \emph{word-length function} for $G$ with respect to $S$.
The left-invariant word-metric $d_{S}$ on $G$ is induced by word-length function, i.e.,
for every $ g, h\in G$,
\[
d_{S}(g, h)=|g^{-1}h|_{S}.
\]

The Cayley graph is the graph whose vertex set is $G$, one vertex for each element in $G$
and any two vertices $ g, h\in G$ are incident with an edge if and only if $g^{-1}h\in S\cup S^{-1}$.

Let $G$ and $N$ be finitely generated groups and let $1_{G}\in G$
and $1_{N}\in N$ be their units. The \emph{support} of a function $f:
N\rightarrow G$ is the set \[
\text{supp}(f)=\{x\in N
|f(x)\neq1_{G}\}.
\]
 The direct sum $\displaystyle\bigoplus_{N}G$ of groups $G$ (or
restricted direct product) is the group of functions
\[
C_{0}(N,G)=\{f: N\rightarrow G \text{ with finite support}\}.
\]
There is a natural action of $N$ on $C_{0}(N,G)$: for all $ a\in N, x\in N, f\in C_{0}(N,G)$,
\[
 a(f)(x) =
f(xa^{-1}).
\]
The semidirect product $C_{0}(N,G)\rtimes N$ is called \emph{restricted wreath
product} and is denoted as $G\wr N$. We recall that the product
in $G\wr N$ is defined by the formula
\[
(f, a)(g, b)=(f a(g), ab) \qquad \forall f,g\in C_{0}(N,G), \forall a,b\in N.
\]
Note that $(f, a)^{-1}=(a^{-1}(f^{-1}),a^{-1}).$

Let $S$ and $T$ be finite generating sets for $G$ and $N$,
respectively. Let $e\in C_{0}(N,G)$ denotes the constant function
with the constant value $1_{G}$. For every $v\in N$ and $b\in G$, let
$
\delta^{b}_{v}:N\rightarrow G
$  be the $\delta$-function, i.e.
\[
\delta^{b}_{v}(v)=b
\text{ and }\delta^{b}_{v}(x)=1_{G} \text{ for }x\neq v.
\]
 Note
that $a(\delta^{b}_{v})=\delta^{b}_{va}$ and hence
$(\delta^{b}_{v} , 1_{N})=(e, v)(\delta^{b}_{1_{N}}, 1_{N})(e, v^{-1})$.
Since every function $f\in C_{0}(N,G)$ can be presented
$\delta^{b_{1}}_{v_{1}}\cdots\delta^{b_{k}}_{v_{k}}$,
\[
(f,
1_{N})=(\delta^{b_{1}}_{v_{1}},1_{N})\cdots(\delta^{b_{k}}_{v_{k}},
1_{N}) \text{ and } (f, u)=(f, 1_{N})(e, u).
\]
 The set
$\widetilde{S}=\{(\delta^{s}_{1_{N}},1_{N}) , (e, t) |s\in S, t\in
T \}$ is a generating set for $G\wr N$.
Note that $G$ and $N$ are subgroups of $G\wr N$.

An explicit formula for the word length of wreath products was found by Parry.
\begin{lem}{\rm(\cite{parry},\cite{SV2007}, Proposition 2.4)}
\label{lem:len}\rm
Let $x=(f,n)\in H \wr \ZZ$, $m=\min\{k\in\ZZ\mid f(k)\neq 1_{H}\}$,
$M=\max\{k\in \ZZ \mid f(k)\neq 1_{H}\}$, then the word-length of $x$
satisfies:
\[
|x|=
\begin{cases}
|n| \qquad &\text{if } f=e, \\
\displaystyle\sum_{i\in \ZZ} |f(i)|+ L_\ZZ(x).\qquad &\text{otherwise.}
\end{cases}
\]
where $e$ is the identity of $\displaystyle\bigoplus_{l\in\ZZ} H$, $L_{\ZZ}(x)$
denotes the length of the shortest path starting from $0$, ending at
$n$ and passing through $m$ and $M$ in the (canonical) Cayley graph
of $\ZZ$.
\end{lem}
\begin{remark}
By the formula of the word-length, we obtain that for every $x=(f,n_{1}),y=(g,n_{2})\in \ZZ \wr \ZZ$,
\[
d(x,y)=|x^{-1}y|=|(n_{1}^{-1}(f^{-1}g),n_{1}^{-1}n_{2})|=\begin{cases}
|n_{1}-n_{2}| \qquad &\text{if } f=g, \\
\displaystyle\sum_{i\in \ZZ} |f(i)-g(i)|+ L_\ZZ(x^{-1}y).\qquad &\text{otherwise.}
\end{cases}
\]
where $L_{\ZZ}(x^{-1}y)$
denotes the length of the shortest path starting from $0$, ending at
$n_{2}-n_{1}$ and passing through all vertices in the support of $n_{1}^{-1}(f^{-1}g)$ in the (canonical) Cayley graph
of $\ZZ$.

Note that $d(x,y)\geq |n_{1}-n_{2}|.$
\end{remark}
\end{subsection}
\begin{subsection}{Saturated Union}
\begin{defi}
Let $\mathcal{U}$ and $\mathcal{V}$ be families of subsets of $X$. The \emph{$r$-saturated union} of $\mathcal{V}$ with $\mathcal{U}$ is defined as
$$\mathcal{V}\cup_r\mathcal{U}=\{N_r(V;\mathcal{U})~|~V\in\mathcal{V}\}\cup\{U\in\mathcal{U}~|~d(U,\mathcal{V})>r\},$$
where $N_r(V;\mathcal{U}) = V \cup\bigcup_{d(U,V)\leq r} U$ and $d(U,\mathcal{V})>r$ means that for every $V\in\mathcal{V}$, $d(U,V)>r$.
\end{defi}

\begin{lem}\rm(see \cite{Bell2011})
\label{r-saturated}
Let $\mathcal{U}$ be an $r$-disjoint and $R$-bounded family of subsets of $X$ with $R\geq r$. Let $\mathcal{V}$ be a $5R$-disjoint, $D$-bounded family of subsets of $X$. Then the family $\mathcal{V}\cup_r\mathcal{U}$ is $r$-disjoint, $(D+2R+2r)$-bounded and $\mathcal{V}\cup_r\mathcal{U}$
covers $\bigcup(\mathcal{V}\cup\mathcal{U})$.

\end{lem}

\end{subsection}

\end{section}

\begin{section}{Main Results}\

\begin{lem}\rm
\label{lem:apc}
Let $a\in\ZZ,k\in\NN$, let $$X=\big(\bigoplus_{\ZZ}\ZZ,[a,a+k]\cap\ZZ\big)=\big\{(f,n)~|~f\in\bigoplus_{\ZZ}\ZZ,n\in\big[a,a+k]\cap\ZZ\}$$
be a subspace of the metric space $\ZZ \wr\ZZ$ with the left-invariant word-metric. For every $m\in\NN$, there exist $B=B(m)>0$ and $B$-bounded families \[\mathcal{U}_{0},\mathcal{U}_{1},\cdots,\mathcal{U}_{(3k+1)2^{3k+1}}\] such that each $\mathcal{U}_i$ is $m$-disjoint for $i\in\{1,2,\cdots, (3k+1)2^{3k+1}\}$, $\mathcal{U}_0$ is $k$-disjoint  and
$\bigcup_{i=0}^{(3k+1)2^{3k+1}}\mathcal{U}_i$ covers $X$.

\end{lem}

\begin{proof}
Without loss of generality , we assume that $m\geq k$.
Let
$$\mathcal{V}_0=\big\{[(2n-1)m,2nm):n\in\mathbb{Z}\big\}\text{ and }\mathcal{V}_1=\big\{[2nm,(2n+1)m):n\in\mathbb{Z}\big\}.$$
Note that $\mathcal{V}_0,\mathcal{V}_1$ are $m$-disjoint and $m$-bounded. Moreover,  $\mathcal{V}_0\cup\mathcal{V}_1$ covers $\mathbb{Z}$.
Let $S=k+m$. For $l=1,...,2^{2m}$, let
$$\mathcal{C}_l=\big\{[(2^{2m}(n-1)+l)2S-m,(2^{2m}n+l)2S-m-k):n\in\mathbb{Z}\big\},$$
$$\mathcal{D}_l=\big\{[(2^{2m}n+l)2S-m-k,(2^{2m}n+l)2S-m):n\in\mathbb{Z}\big\}.$$
Note that each $\mathcal{C}_l$ is $k$-disjoint and $2^{3m+3}$-bounded, $\mathcal{D}_l$ is $m$-disjoint and $k$-bounded, $\mathcal{D}_{l}\cup\mathcal{C}_{l}$ covers $\mathbb{Z}$.
Moreover, $\bigcup_{l=1}^{2^{2m}}\mathcal{D}_l$ is $m$-disjoint. Let

$$\mathcal{W}_l=\big\{(V_1,V_2,\cdots,V_{2m})~|~V_i\in \mathcal{V}_{\phi(l)_i},i\in\{1,...,2m\}\big\}$$
where $\phi$ is a bijection from $\{1,...,2^{2m}\}$ to $\{0,1\}^{2m}$.
Then  \[\text{$\bigcup_{l=1}^{2^{2m}}\mathcal{W}_l$ is disjoint. i.e., for every $W_{1},W_{2}\in\bigcup_{l=1}^{2^{2m}}\mathcal{W}_l$ and $W_{1}\neq W_{2}$,
$W_{1}\cap W_{2}=\varnothing.$}\]

Let
\begin{multline*}
\mathcal{U}_{0}=\{(\{\widetilde{x}\}\times W_{1}\times\prod_{i=a-k}^{a+2k}C_i\times W_{2}\times \{\widetilde{y}\},[a,a+k]\cap\mathbb{Z})~|\\
\big(\prod_{i=a-k}^{a+2k}C_i,(W_{1},W_{2})\big)\in\bigcup_{l=1}^{2^{2m}}(\mathcal{C}_l^{3k+1}\times \mathcal{W}_l),\widetilde{x},\widetilde{y}\in \bigoplus_{\ZZ}\ZZ.\},
\end{multline*}
\begin{multline*}
\mathcal{U}_{2^{3k+1}(s-a+k)+t}=\{(\{\widetilde{x}\}\times W_{1}\times\prod_{i=a-k}^{s-1}V_i\times D_s\times \prod_{i=s+1}^{a+2k}V_i\times W_{2}\times \{\widetilde{y}\},[a,a+k]\cap\mathbb{Z})~| \\V_{i}\in \mathcal{V}_{\rho(t)_i},i\in\{a-k,a-k+1,...,s-1,s+1,...,a+2k\},
\big(D_s,(W_{1},W_{2})\big)\in\bigcup_{l=1}^{2^{2m}}(\mathcal{D}_l\times \mathcal{W}_l),\widetilde{x},\widetilde{y}\in \bigoplus_{\ZZ}\ZZ.\}
\end{multline*}
where $\rho$ is a bijection from $\{1,...,2^{3k+1}\}$ to $\{0,1\}^{3k+1}$, $s\in\{a-k,a-k+1,...,a+2k\}$ and $t\in\{1,2,3,...,2^{3k+1}\}$.

\begin{itemize}
\item First we will prove that $\mathcal{U}_0$ is $k$-disjoint and $B(m)$-bounded, where $B(m)=12m+2m^{2}+2^{6m+4}$.

Let
\[U=\big(\{\widetilde{x}\}\times W_{1}\times\prod_{i=a-k}^{a+2k}C_i\times W_{2}\times \{\widetilde{y}\},[a,a+k]\cap\mathbb{Z}\big)
\]
and
\[U^{'}=\big(\{\widetilde{x^{'}}\}\times W_{1}^{'}\times\prod_{i=a-k}^{a+2k}C_i^{'}\times W_{2}^{'}\times \{\widetilde{y^{'}}\},[a,a+k]\cap\mathbb{Z}\big).
\]
Assume that $U,U^{'}\in\mathcal{U}_0$ and $U\neq U^{'}$. For every $x=(f,n_{1})\in U,y=(g,n_{2})\in U^{'}$,
\begin{itemize}
\item Case 1. $(\{\widetilde{x}\}, W_{1}, W_{2},\{\widetilde{y}\})\neq(\{\widetilde{x^{'}}\}, W_{1}^{'}, W_{2}^{'},\{\widetilde{y^{'}}\}).$

Then $\bigcup_{l=1}^{2^{2m}}\mathcal{W}_l$ is disjoint implies that there exists $j<a-k$ or $j>a+2k$ such that $f(j)\neq g(j)$. i.e., $(f^{-1}g)(j)\neq 0$, which implies that
there exists $j<a-k-n_{1}\leq -k$ or $j>a+2k-n_{1}\geq k$ such that $n_{1}^{-1}(f^{-1}g)(j)\neq 0$.
It follows that $L_{\ZZ}(x^{-1}y)\geq k$. And hence $d(x,y)\geq k$.

\item Case 2. $(\{\widetilde{x}\}, W_{1}, W_{2},\{\widetilde{y}\})=(\{\widetilde{x^{'}}\}, W_{1}^{'}, W_{2}^{'},\{\widetilde{y^{'}}\}).$

Since $\bigcup_{l=1}^{2^{2m}}\mathcal{W}_l$ is disjoint, $(W_{1}, W_{2})= (W_{1}^{'}, W_{2}^{'})\in \mathcal{W}_l$ for some unique $l\in\{1,2,\cdots,2^{2m}\}$. Then
\[
\prod_{i=a-k}^{a+2k}C_i, \prod_{i=a-k}^{a+2k}C_i^{'}\in \mathcal{C}_l^{3k+1}\text{~and~}\prod_{i=a-k}^{a+2k}C_i\neq\prod_{i=a-k}^{a+2k}C_i^{'}
\]
It follows that there exists $i\in [a-k,a+2k]\cap\ZZ$ such that $C_i\neq C_i^{'}$ and $C_i,C_i^{'}\in \mathcal{C}_l$, which implies that
$|f(i)-g(i)|\geq k.$  And hence $d(x,y)\geq k$.

\end{itemize}
So $d(U,U^{'})\geq k$. So $\mathcal{U}_0$ is $k$-disjoint.

For every $x=(f,n_{1}),z=(h,n_{3})\in U$.
\begin{itemize}
\item If $f=h$, then $d(x,z)=|n_{3}-n_{1}|\leq k\leq m$.
\item  Otherwise, since
\[
\displaystyle\sum_{i\in \ZZ} |f(i)-g(i)|\leq 2m^{2}+(3k+1)2^{3m+3}\leq 2m^{2}+(3m+1)2^{3m+3}\leq 2m^{2}+2^{6m+4},
\]
\[
\text{~and~}\text{supp~}n_{1}^{-1}(f^{-1}g)\subseteq[a-k-m-n_{1},a+2k+m-n_{1}]\subseteq[-2k-m,2k+m]\subseteq[-3m,3m]
\]
implies that $L_{\ZZ}(x^{-1}z)\leq 12m,$ $d(x,z)\leq 12m+2m^{2}+2^{6m+4}=B(m)$.
\end{itemize}
So diam~$U\leq B(m), \forall~ U\in \mathcal{U}_0$, i.e., $\mathcal{U}_0$ is $B(m)$-bounded.

\item Now we will prove that for every $i\in\{1,2,\cdots,(3k+1)2^{3k+1}\}$, $\mathcal{U}_i$ is $m$-disjoint and $B(m)$-bounded.

Let
\[U=(\{\widetilde{x}\}\times W_{1}\times\prod_{i=a-k}^{s-1}V_i\times D_s\times \prod_{i=s+1}^{a+2k}V_i\times W_{2}\times \{\widetilde{y}\},[a,a+k]\cap\mathbb{Z})
\]
and
\[U^{'}=(\{\widetilde{x^{'}}\}\times W_{1}^{'}\times\prod_{i=a-k}^{s-1}V_i^{'}\times D_s^{'}\times \prod_{i=s+1}^{a+2k}V_i^{'}\times W_{2}^{'}\times \{\widetilde{y}\},[a,a+k]\cap\mathbb{Z}).
\]
Assume that $U,U^{'}\in\mathcal{U}_i$ and $U\neq U^{'}$. For every $x=(f,n_{1})\in U,y=(g,n_{2})\in U^{'}$,
\begin{itemize}
\item case 1. $ D_s\neq D_s^{'}.$

Since $\bigcup_{l=1}^{2^{2m}}\mathcal{D}_l$ is $m$-disjoint, $|f(s)-g(s)|\geq m$.
So $d(x,y)\geq m$.

\item case 2. $ D_s=D_s^{'}.$

Then $ D_s=D_s^{'}\in \mathcal{D}_l$ for some unique $l\in\{1,2,\cdots,2^{2m}\}$. By definition of~ $\mathcal{U}_i$,
$(W_{1}, W_{2}),(W_{1}^{'}, W_{2}^{'})\in \mathcal{W}_l$.
If $(W_{1}, W_{2})\neq(W_{1}^{'}, W_{2}^{'})$, then $\displaystyle\sum_{i\in \ZZ} |f(i)-g(i)|\geq m.$
If $(W_{1}, W_{2})=(W_{1}^{'}, W_{2}^{'})$, then
\[
(\{\widetilde{x}\},\prod_{i=a-k}^{s-1}V_i, \prod_{i=s+1}^{a+2k}V_i,\{\widetilde{y}\})\neq(\{\widetilde{x^{'}}\},\prod_{i=a-k}^{s-1}V_i^{'}, \prod_{i=s+1}^{a+2k}V_i^{'},\{\widetilde{y^{'}}\}).
\]
It follows that at least one of the following two situations holds.
\begin{itemize}
\item [(a)] There exists $j\in\{a-k,a-k+1,\cdots,a+2k\}$ such that $V_{j}\neq V_{j}^{'}\in \mathcal{V}_{\rho(t)_j}$.
Then $|f(j)-g(j)|\geq m$.
So $d(x,y)\geq m$.
\item [(b)] There exists $j<a-k-m$ or $j>a+2k+m$ such that $f(j)\neq g(j)$. i.e., $(f^{-1}g)(j)\neq 0$, which implies that
there exists $j<a-k-m-n_{1}\leq -m$ or $j>a+2k+m-n_{1}\geq m$ such that $n_{1}^{-1}(f^{-1}g)(j)\neq 0$.
It follows that $L_{\ZZ}(x^{-1}y)\geq m$. So $d(x,y)\geq m$.
\end{itemize}

\end{itemize}
Then $d(U,U^{'})\geq m$. So $\mathcal{U}_i$ is $m$-disjoint

For every $x=(f,n_{1}),z=(h,n_{3})\in U$.
\begin{itemize}
\item If $f=h$, then $d(x,z)=|n_{3}-n_{1}|\leq k\leq m$.
\item  Otherwise, since
\[
\displaystyle\sum_{i\in \ZZ} |f(i)-g(i)|\leq 2m^{2}+(3k+1)m\leq 2m^{2}+(3m+1)m.
\]
\[
\text{~and~}\text{supp~}n_{1}^{-1}(f^{-1}g)\subseteq[a-k-m-n_{1},a+2k+m-n_{1}]\subseteq[-2k-m,2k+m]\subseteq[-3m,3m]
\]
implies that $L_{\ZZ}(x^{-1}z)\leq 12m,$ $d(x,z)\leq 2m^{2}+(3m+1)m+12m\leq B(m).$
\end{itemize}
So diam~$U\leq B(m).$ i.e., $\mathcal{U}_i$ is $B(m)$-bounded.

\item Finally, we will prove that $\bigcup_{i=0}^{(3k+1)2^{3k+1}}\mathcal{U}_i$ covers $X$.

Indeed, let $x=(f,n)\in X\setminus \bigcup\mathcal{U}_0$,
then there exists unique $l\in \{1,...,2^{2m}\}$ and $(W_{1},W_{2})\in \mathcal{W}_l$ such that
 \[\text{$\big(f(a-k-m),\cdots,f(a-k-1)\big)\in W_{1}$ and
$\big(f(a+2k+1),\cdots,f(a+2k+m)\big)\in W_{2}$.}\] Since $x\notin \bigcup\mathcal{U}_0$, we have $\big(f(a-k),f(a-k+1)\cdots,f(a+2k)\big)\notin \bigcup\{\prod_{i=a-k}^{a+2k}C_i~|~C_i\in\mathcal{C}_l\}$. Then there exists $s\in\{a-k,a-k+1,...,a+2k\}$ such that $x_{s}\notin \bigcup\mathcal{C}_l$. Since $\mathcal{C}_l\bigcup\mathcal{D}_l$ covers $\mathbb{Z}$, there exists $D_s\in \mathcal{D}_l$ such that $x_{s}\in D_s$. Since $\mathcal{V}_0\bigcup\mathcal{V}_1$ covers $\mathbb{Z}$, we may take $t\in\{1,...,2^{3k+1}\}$ such that $\big(f(a-k),f(a-k+1)\cdots,f(a+2k)\big)\in \prod_{i=a-k}^{a+2k}V_i$, where $V_{i}\in\mathcal{V}_{\rho(t)_i},i\in\{a-k,a-k+1,...,a+2k\}$. So
$x\in\bigcup\mathcal{U}_{2^{3k+1}(s-a+k)+t}.$

\end{itemize}

\end{proof}

\begin{thm}\rm
Let $X$ be the metric space $\ZZ \wr\ZZ$ with the left-invariant word-metric. Then trasdim$X\leq \omega+1.$
Consequently, $\ZZ \wr\ZZ$ has asymptotic property $C$.

\end{thm}

\begin{proof}
By Lemma \ref{lem:trasdim}, it suffices to show that for every $k,m\in\NN$,  there are uniformly bounded families $\mathcal{V}_{-1},\mathcal{V}_{0},\cdots,\mathcal{V}_{(6k+2)2^{3k+1}}$ such that $\mathcal{V}_i$ is $k$-disjoint for $i=-1,0$, $\mathcal{V}_j$ is $m$-disjoint for $j=1,2,\cdots, (6k+2)2^{3k+1}$ and
$\bigcup_{i=-1}^{(6k+2)2^{3k+1}}\mathcal{V}_i$ covers $X$.

Without loss of generality , we assume that $m= kp$ for some $p\in\NN$.
For $i\in\ZZ$, let
$
A_{i}=[ik,ik+k]\cap\ZZ.
$
By Lemma \ref{lem:apc}, there exist $B_{1}(m)>0$ and $B_{1}(m)$-bounded families
\[\mathcal{U}_{0}^{ip+1},\mathcal{U}_{1}^{ip+1},\cdots,\mathcal{U}_{(3k+1)2^{3k+1}}^{ip+1}\] satisfying each $\mathcal{U}_j^{ip+1}$ is $m$-disjoint for $j\in\{1,2,\cdots, (3k+1)2^{3k+1}\}$, $\mathcal{U}_0^{ip+1}$ is $k$-disjoint  and
$\bigcup_{j=0}^{(3k+1)2^{3k+1}}\mathcal{U}_j^{ip+1}$ covers $\big(\bigoplus_{\ZZ}\ZZ,A_{ip+1}\big)$.
Similarly, there exist $D_{1}(m)>0$ and $D_{1}(m)$-bounded families
\[\mathcal{U}_{0}^{ip+2},\mathcal{U}_{1}^{ip+2},\cdots,\mathcal{U}_{(3k+1)2^{3k+1}}^{ip+2}\] satisfying each $\mathcal{U}_j^{ip+2}$ is $5B_{1}(m)$-disjoint for $j\in\{1,2,\cdots, (3k+1)2^{3k+1}\}$, $\mathcal{U}_0^{ip+2}$ is $k$-disjoint  and
$\bigcup_{j=0}^{(3k+1)2^{3k+1}}\mathcal{U}_j^{ip+2}$ covers $\big(\bigoplus_{\ZZ}\ZZ,A_{ip+2}\big)$.
For $j\in\{1,2,\cdots, (3k+1)2^{3k+1}\}$, let
\[
\mathcal{V}_j^{ip+3}=\mathcal{U}_j^{ip+2}\cup_m\mathcal{U}_j^{ip+1}.
\]
Then $\mathcal{V}_j^{ip+3}$ is $m$-disjoint and $B_{2}(m)$-bounded by Lemma \ref{r-saturated}.
Similarly, there exist $D_{2}(m)>0$ and $D_{2}(m)$-bounded families
\[\mathcal{U}_{0}^{ip+3},\mathcal{U}_{1}^{ip+3},\cdots,\mathcal{U}_{(3k+1)2^{3k+1}}^{ip+3}\] satisfying each $\mathcal{U}_j^{ip+3}$ is $5B_{2}(m)$-disjoint for $j\in\{1,2,\cdots, (3k+1)2^{3k+1}\}$, $\mathcal{U}_0^{ip+3}$ is $k$-disjoint  and
$\bigcup_{j=0}^{(3k+1)2^{3k+1}}\mathcal{U}_j^{ip+3}$ covers $\big(\bigoplus_{\ZZ}\ZZ,A_{ip+3}\big)$.
For $j\in\{1,2,\cdots, (3k+1)2^{3k+1}\}$, let
\[
\mathcal{V}_j^{ip+4}=\mathcal{U}_j^{ip+3}\cup_m\mathcal{V}_j^{ip+3}.
\]
Then $\mathcal{V}_j^{ip+4}$ is $m$-disjoint and $B_{3}(m)$-bounded by Lemma \ref{r-saturated}.
After finite steps, we obtain that $m$-disjoint and $B_{p-1}(m)$-bounded families
\[
\mathcal{V}_1^{ip+p},\mathcal{V}_2^{ip+p},\cdots,\mathcal{V}_{(3k+1)2^{3k+1}}^{ip+p}.
\]
Note that $\mathcal{V}_j^{ip+p}$ is a family of subsets of $\big(\bigoplus_{\ZZ}\ZZ,\bigcup_{j=1}^{p}A_{ip+j}\big)$ and
\[
\big(\bigcup_{j=1}^{p}\mathcal{U}_{0}^{ip+j}\big)\bigcup\big(\bigcup_{j=1}^{(3k+1)2^{3k+1}}\mathcal{V}_{j}^{ip+p}\big)\text{ covers }\big(\bigoplus_{\ZZ}\ZZ,\bigcup_{j=1}^{p}A_{ip+j}\big)
\]
For $j\in\{1,2,\cdots, (3k+1)2^{3k+1}\}$, let
\[
\mathcal{V}_j=\bigcup_{n\in\ZZ}\mathcal{V}_j^{2np+p}\text{ and }\mathcal{V}_{j+(3k+1)2^{3k+1}}=\bigcup_{n\in\ZZ}\mathcal{V}_j^{(2n+1)p+p}.
\]

Since for every $i_{1},i_{2}\in\ZZ$ and $i_{1}\neq i_{2},$
\[
d\big(\big(\bigoplus_{\ZZ}\ZZ,\bigcup_{j=1}^{p}A_{2i_{1}p+j}\big),\big(\bigoplus_{\ZZ}\ZZ,\bigcup_{j=1}^{p}A_{2i_{2}p+j}\big)\big)\geq m \text{~and~}d\big((\bigoplus_{\ZZ}\ZZ,\bigcup_{j=1}^{p}A_{(2i_{1}+1)p+j}),(\bigoplus_{\ZZ}\ZZ,\bigcup_{j=1}^{p}A_{(2i_{2}+1)p+j})\big)\geq m,
\]
each $\mathcal{V}_j$ is $m$-disjoint and uniformly bounded for $j\in\{1,2,\cdots, (6k+2)2^{3k+1}\}$.
Let
\[
\mathcal{V}_0=\bigcup_{n\in\ZZ}\mathcal{U}_0^{2n}\text{ and }\mathcal{V}_{-1}=\bigcup_{n\in\ZZ}\mathcal{U}_0^{2n+1}.
\]
Since
\[
\{\big(\bigoplus_{\ZZ}\ZZ,A_{2n}\big)~|~n\in\ZZ\} \text{~is }k\text{-disjoint and } \mathcal{U}_0^{2n} \text{ is a family of subsets of } \big(\bigoplus_{\ZZ}\ZZ,A_{2n}\big) \text{~and~}\]
\[\{\big(\bigoplus_{\ZZ}\ZZ,A_{2n+1}\big)~|~n\in\ZZ\} \text{~is }k\text{-disjoint and } \mathcal{U}_0^{2n+1} \text{ is a family of subsets of } \big(\bigoplus_{\ZZ}\ZZ,A_{2n+1}\big),
\]
 $\mathcal{V}_0$ and $\mathcal{V}_{-1}$ is $k$-disjoint and uniformly bounded.
 
Finally, we will prove that $\bigcup_{i=-1}^{(6k+2)2^{3k+1}}\mathcal{V}_i$ covers $\ZZ \wr\ZZ$.

Since 
\[
\ZZ \wr\ZZ=\big(\bigoplus_{\ZZ}\ZZ,\ZZ\big)=\bigcup_{i\in\ZZ}\big(\bigoplus_{\ZZ}\ZZ,A_{i}\big)=\bigcup_{i\in\ZZ}\bigcup_{j=1}^{p}(\bigoplus_{\ZZ}\ZZ,A_{ip+j}\big)
=\bigcup_{i\in\ZZ}(\bigoplus_{\ZZ}\ZZ,\bigcup_{j=1}^{p}A_{ip+j}\big)
\]
and
\[
\big(\bigcup_{j=1}^{p}\mathcal{U}_{0}^{ip+j}\big)\bigcup\big(\bigcup_{j=1}^{(3k+1)2^{3k+1}}\mathcal{V}_{j}^{ip+p}\big)\text{ covers }\big(\bigoplus_{\ZZ}\ZZ,\bigcup_{j=1}^{p}A_{ip+j}\big),
\]
we obtain that
\[
\big(\bigcup_{i\in\ZZ}\bigcup_{j=1}^{p}\mathcal{U}_{0}^{ip+j}\big)\bigcup\big(\bigcup_{i\in\ZZ}\bigcup_{j=1}^{(3k+1)2^{3k+1}}\mathcal{V}_{j}^{ip+p}\big)\text{ covers }\ZZ \wr\ZZ.
\]
Note that 
\[
\bigcup_{i\in\ZZ}\bigcup_{j=1}^{p}\mathcal{U}_{0}^{ip+j}=\bigcup_{i\in\ZZ}\mathcal{U}_{0}^{i}=\big(\bigcup_{n\in\ZZ}\mathcal{U}_0^{2n}\big)\bigcup\big(\bigcup_{n\in\ZZ}\mathcal{U}_0^{2n+1}\big)
=\mathcal{V}_0\cup\mathcal{V}_{-1}
\]
and
\[
\begin{split}
\bigcup_{i\in\ZZ}\bigcup_{j=1}^{(3k+1)2^{3k+1}}\mathcal{V}_{j}^{ip+p}&=\big(\bigcup_{j=1}^{(3k+1)2^{3k+1}}\bigcup_{n\in\ZZ}\mathcal{V}_{j}^{2np+p}\big)\bigcup
\big(\bigcup_{j=1}^{(3k+1)2^{3k+1}}\bigcup_{n\in\ZZ}\mathcal{V}_{j}^{(2n+1)p+p}\big)\\
&=\big(\bigcup_{j=1}^{(3k+1)2^{3k+1}}\mathcal{V}_{j}\big)\bigcup
\big(\bigcup_{j=1}^{(3k+1)2^{3k+1}}\mathcal{V}_{j+(3k+1)2^{3k+1}}\big)\\
&=\bigcup_{j=1}^{(6k+2)2^{3k+1}}\mathcal{V}_{j}.
\end{split}
\]
Therefore, $\bigcup_{i=-1}^{(6k+2)2^{3k+1}}\mathcal{V}_i$ covers $\ZZ \wr\ZZ$.

\end{proof}

\end{section}

{\bf Acknowledgments.} The authors wish to thank the reviewers for careful
reading and valuable comments. This work was supported by NSFC grant of P.R. China (No.12071183,11871342).


\providecommand{\bysame}{\leavevmode\hbox to3em{\hrulefill}\thinspace}
\providecommand{\MR}{\relax\ifhmode\unskip\space\fi MR }
\providecommand{\MRhref}[2]{%
  \href{http://www.ams.org/mathscinet-getitem?mr=#1}{#2}
}
\providecommand{\href}[2]{#2}

\end{document}